\documentclass[10pt]{amsart}
\usepackage{palatino, amssymb, amsmath, amsfonts, latexsym, amscd}
\input xy

\catcode`\@=11

\newlength{\tabwidth}
\newlength{\tabheight}
\newlength{\tabrule}
\newlength{\tabwidthx}
\newlength{\tabheightx}

\def\gentabbox#1#2#3#4{\vbox to \tabheight{\setlength{\tabrule}{#3}%
  \setlength{\tabwidthx}{#1\tabwidth}\addtolength{\tabwidthx}{\tabrule}%
  \setlength{\tabheightx}{#2\tabheight}\addtolength{\tabheightx}{-\tabheight}%
  \hbox to #1\tabwidth{%
    \hspace{-0.5\tabrule}\rule{\tabrule}{#2\tabheight}\hspace{-\tabrule}%
    \vbox to #2\tabheight{\hsize=\tabwidthx%
      \vspace{-0.5\tabrule}\hrule width\tabwidthx height\tabrule%
      \vspace{-0.5\tabrule}\vss%
      \hbox to \tabwidthx{\hss#4\hss}%
        \vss\vspace{-0.5\tabrule}%
      \hrule width\tabwidthx height\tabrule\vspace{-0.5\tabrule}}%
    \hspace{-\tabrule}\rule{\tabrule}{#2\tabheight}\hspace{-0.5\tabrule}}%
  \vspace{-\tabheightx}}}
\def\genblankbox#1#2{\vbox to \tabheight{\vfil\hbox to #1\tabwidth{\hfil}}}
\def\tabbox#1#2#3{\gentabbox{#1}{#2}{0.4pt}{\strut #3}}

\catcode`\:=13 \catcode`\.=13 \catcode`\;=13 \catcode`\>=13
\catcode`\^=13
\def:#1\\{\hbox{$#1$}}
\def.#1{\tabbox{1}{1}{$#1$}}
\def>#1{\tabbox{2}{1}{$#1$}}
\def^#1{\tabbox{1}{2}{$#1$}}
\def;{\genblankbox{1}{1}\relax}
\catcode`\:=12 \catcode`\.=12 \catcode`\;=12 \catcode`\>=12
\catcode`\^=7

\newenvironment{tableau}{\bgroup\catcode`\:=13 \catcode`\.=13
  \catcode`\;=13 \catcode`\>=13 \catcode`\^=13
  \def\b##1##2##3{\gentabbox{##1}{##2}{1.2pt}{\vbox{##3}}}%
  \def\n##1##2##3{\gentabbox{##1}{##2}{0.4pt}{\vbox{##3}}}%
  \vcenter\bgroup\offinterlineskip}{\egroup\egroup}

\catcode`\@=11

\xyoption{all}

\newtheorem{theorem}{Theorem}[section]
\newtheorem{lemma}[theorem]{Lemma}
\newtheorem{prop}[theorem]{Proposition}
\newtheorem{cor}[theorem]{Corollary}

\theoremstyle{definition}
\newtheorem{definition}[theorem]{Definition}
\newtheorem{example}[theorem]{Example}

\theoremstyle{remark}
\newtheorem{remark}[theorem]{Remark}
\newcommand{\U}{\dot{\mathbf U}}
\newcommand{\B}{\dot{\mathbf B}}
\newcommand{\ca}{\{A\}}
\newcommand{\cb}{\{B\}}
\newcommand{\cc}{\{C\}}
\newcommand{\seq}{\mathfrak S_{D,n}}
\newcommand{\mat}{\mathfrak S_{D,n,n}}
\numberwithin{equation}{section}

\begin{document}

\title{Cells in quantum affine $\mathfrak{sl}_n$}

\author{Kevin McGerty}
\address{Department of Mathematics, Massachusetts Institute of Technology. }

\date{September 4, 2002}

\begin{abstract}
We study Lusztig's theory of cells for quantum affine $\mathfrak{sl}_n$. Using the
geometric construction of the quantum group due to Lusztig and Ginzburg--Vasserot, we
describe explicitly the two-sided cells, the number of left cells in a two--sided cell,
and the asymptotic algebra, verifying conjectures of Lusztig.
\end{abstract}

\maketitle

\section{Introduction}

Given any algebra with a specified basis it is
possible to define a notion of cell from the set of ideals, left, right or two-sided,
which are spanned by a subset of the basis. Of course, if one picks the basis arbitrarily,
 it is unlikely that these objects will contain any interesting information about the algebra
 in question. However, if the algebra has a natural choice of basis, the situation can be quite
 different. Examples of this arise in a number of places: The cells
attached to the Kazhdan--Lusztig basis of a Hecke algebra associated to a finite Weyl group turn
out to be crucial in the classification of the characters of finite groups of Lie type. On the other hand, although
 the plus part of the quantum group possesses a natural ``canonical basis'',
 the theory of cells there is trivial. If we extend the canonical basis
 to one for the modified quantum group $\dot{\mathbf U}$ however, the
 theory of cells is once again interesting.

In the case of quantum groups of finite type, work of Lusztig
\cite{L95} completely describes the cells. In that paper Lusztig also gave a
conjectural description of the cell structure in the case of (degenerate) affine
quantum groups. In this paper we show that the geometric construction of
 $\dot{\mathbf U}$ in  \cite{L99} can be used to give complete information about
 the cell structure of quantum affine $\mathfrak{sl}_n$.
 Just as in the case of affine Weyl groups, the cells are closely
 related to the finite dimensional representation theory of the algebra.
 We will first investigate the structure of cells in the ``affine q-Schur
 algebra'' $\mathfrak A_D$ and then show how this can be used to obtain
 the cell structure of $\U$.

We begin by recalling the definition of cells. Suppose $R$ is a
ring, and $A$ an associative algebra over $R$, with an $R$--basis
$B$. We say that a left ideal is \textit{based} if it is the span of a
subset of the basis $B$. We define a preorder on the elements of
$B$ as follows. Let $x \preceq_L y$ for $x,y \in B$ if $x$ lies in
every based left ideal which contains $y$. The equivalence classes of
this preorder are precisely the left cells of $A$. If we replace ``left ideal''
with ``right ideal'' or ``two--sided ideal'' we get the corresponding
notion of right cells or two--sided cells.

\section{The Affine q-Schur Algebra} \label{qschur}

We first give a description of the affine q-Schur algebra
 $\mathfrak A_D$ as a commutator algebra, following the notation of
 \cite{L99} and \cite{M02} (see also \cite{ScV}, \cite{GV93}). Thus let $V$ be a free rank $D$ module over
 $\mathbf k[\epsilon,\epsilon^{-1}]$, where $\mathbf k$ is a finite field
 of $q$ elements, and $\epsilon$ is an indeterminate.

Let $\mathcal F^n$ be the space of n-step periodic lattices, i.e. sequences $\mathbf L = (L_i)_{i \in \mathbb Z}$
 of lattices in our free module $V$ such that $L_i \subset L_{i+1}$, and $L_{i-n}=\epsilon L_i$.
 The group $G = \mathrm{Aut}(V)$ acts on $\mathcal F^n$ in the natural way. Let $\seq$ be the set of nonnegative integer
 sequences $(a_i)_{i \in \mathbb Z}$, such that $a_i = a_{i+n}$ and $\sum_{i=1}^n a_i = D$, and let $\mat$ be
 the set of $\mathbb Z \times \mathbb Z$ matrices $A = (a_{i,j})_{i,j \in \mathbb Z}$
 with nonnegative entries such that $a_{i,j}= a_{i+n,j+n}$ and $\sum_{i \in [1,n], j \in \mathbb Z} a_{i,j} = D$.
 The orbits of $G$ on $\mathcal F^n$ are indexed by $\seq$, where $\mathbf L$ is in the orbit $\mathcal F_{\mathbf a}$
 corresponding to $\mathbf a$ if $a_i = \mathrm{dim}_{\mathbf k}(L_i/L_{i-1})$. The orbits of $G$ on
 $\mathcal F^n \times \mathcal F^n$ are indexed by the matrices $\mat$, where a pair $(\mathbf L, \mathbf L')$
 is in the orbit $\mathcal O_A$ corresponding to $A$ if
 \[
 a_{i,j}=\mathrm{dim}\biggl(\frac{L_i\cap L_j'}{(L_{i-1}\cap
L_j')+(L_i\cap L_{j-1}')}\biggr).
 \]
For $A \in \mat$ let $r(A), c(A) \in \seq$ be given by $r(A)_i = \sum_{j \in \mathbb Z} a_{i,j}$ and
 $r(A)_j = \sum_{i \in \mathbb Z} a_{i,j}$.

Similarly let $\mathcal B^D$ be the space of complete periodic lattices, that is, sequences of lattices
 $\mathbf L = (L_i)$ such that $L_i \subset L_{i+1}$, $L_{i-D}=\epsilon L_i$, and
 $\rm{dim}_{\mathbf k}(L_i/L_{i-1})=1$ for all $i \in\mathbb Z$.  Let $\mathbf b_0 = (\ldots,1,1,\ldots)$.
 The orbits of $G$ on $\mathcal B^D\times \mathcal B^D$ are indexed by matrices $A \in \mathfrak S_{n,n,n}$
 where the matrix $A$ must have $r(A)=c(A)=\mathbf b_0$.

Let $\mathfrak A_{D,q}$, $\mathfrak H_{D,q}$ and $\mathfrak T_{D,q}$ be the span of the characteristic functions
 of the $G$ orbits on $\mathcal F^n \times \mathcal F^n$, $\mathcal B^D \times \mathcal B^D$ and $\mathcal F^n
 \times \mathcal B^D$ respectively. Convolution makes $\mathfrak A_{D,q}$ and $\mathfrak H_{D,q}$ into algebras and
 $\mathfrak T_D$ into a $\mathfrak A_{D,q}$--$\mathfrak H_{D,q}$ bimodule.
 For $A \in \mat$ set
 \[
 d_A = \sum_{i \geq k, j<l, 1\leq i \leq n}a_{ij}a_{kl}.
 \]
 Let $\{[A] \colon A \in \mat\}$ be the basis of $\mathfrak A_{D,q}$ given by $q^{-d_A/2}$ times the characteristic
 function of the orbit corresponding to $A$. When $D=n$ an obvious subset of this basis spans $\mathfrak H_{D,q}$.

All of these spaces of functions are the specialization at $v = \sqrt q$ of modules over $\mathcal A =
 \mathbb Z[v,v^{-1}]$, which we denote by $\mathfrak A_D$, $\mathfrak H_D$ and $\mathfrak T_D$ respectively
 (here $v$ is an indeterminate). The $\mathcal A$--algebra $\mathfrak A_D$ is the ``affine q-Schur algebra'',
 it is easy to see that it is the commutator algebra of the right $\mathcal H_D$--module $\mathfrak T_D$.
 This is the affine version of the geometric Schur--Weyl duality first described in \cite{GL93}.
Recall from \cite[section 4]{L99} that $\mathfrak A_{D}$ possesses
a canonical basis $\mathfrak B_D$ consisting of elements $\{A\}$, for
$A\in \mathfrak S_{D,n,n}$. We have
\[
\{A\} = \sum_{A_1; A_1 \leq A}\Pi_{A_1,A}[A_1].
\]
where $\leq $ is a natural partial order on $\mat$ and the $\Pi_{A_1,A}$ are certain
polynomials in $\mathbb Z[v^{-1}]$.
$\mathfrak A_D$ has a natural antiautomorphism $\Psi$ which sends $[A]$ to $[A^t]$, and a related
antiautomorphism $\rho$ such that $\rho([A]) = v^{d_A - d_{A^t}}[A^t]$.
In \cite{M02} a natural inner product $(\cdot,\cdot)_D$  is defined on $\mathfrak A_D$ (this is different
 from the inner product given in \cite{L99}). It has the property that
\[
(xy,z)_D = (y,\rho(x)z)_D; \qquad x,y,z \in \mathfrak A_D.
\]
Moreover it can be shown that
\[
(\{A\},\{B\})_D \in \delta_{A,B} + v^{-1}\mathbb N[v^{-1}].
\]
so that the inner product is positive with respect to the canonical basis, and the
canonical basis is ``almost orthogonal''.
\bigskip

Let $\mathcal H_D$ be the affine Hecke algebra of $GL_D$, thus $\mathcal H_D$ is an algebra over
 $\mathbb Z[v,v^{-1}]$ generated by symbols $T_i, X_j,X_j^{-1}$, where $i \in \{1,2,\ldots ,D-1\}$,
 and $j \in \{1,2,\ldots,D\}$, subject to the relations

\begin{itemize}
\item $(T_i -v)(T_i+v)=0, \qquad T_iT_{i+1}T_i =
T_{i+1}T_iT_{i+1}$, for $i = 1,2,\ldots,D-1$;

\item $T_iT_j=T_jT_i$ if $|i-j| \geq 2$;

\item $X_iX_i^{-1}=X_i^{-1}X_i=1, \qquad X_iX_j=X_jX_i$, for all
$i,j$;

\item $T_i^{-1}X_iT_i^{-1}=X_{i+1}$ for $ i=1,2,\ldots,D-1; \qquad
T_iX_j=X_jT_i$ for $j \neq i,i+1$.

\end{itemize}
This is the ``Bernstein presentation''. Let $W$ be the affine Weyl
group of type $GL_D$ (when we wish to specify D, we will use the notation
$\hat{A}_{D-1}$), that is, $W$ is the semidirect product of
the symmetric group $S_D$ with $\mathbb Z^D$. It is an extension
by $\mathbb Z$ of a Coxeter group, thus the usual yoga can be
used to extend the Kazhdan-Lusztig theory. Let the set of simple
reflections of the Coxeter group be $S = \{s_i\colon i = 0,1,\dots,D\}$.
The ``Iwahori presentation'' of $\mathcal H_D$ yields a basis $\{T_w\colon w \in W\}$.
Lusztig observed that $W$ has a natural incarnation as a
permutation group on the integers, indeed $W$ is isomorphic to the set of
all permutations $\sigma$ of the integers such that $\sigma(i+D)=\sigma(i)+D$.
See for example \cite{Xi} for more details. Thus an element
of $W$ obviously corresponds to an infinite permutation matrix,
which we denote $A_w$ when we wish to make the distinction between
the group element and the matrix. These permutation matrices are precisely the matrices
indexing the $G$-orbits on $\mathcal B^D \times \mathcal B^D$ described above.

\begin{prop}
The map $\mathcal H_{v=\sqrt q} \to \mathfrak H_{D,q}$ which sends
$T_w \mapsto [A_w]$ is an algebra isomorphism. \hfill $\square$
\end{prop}

\bigskip

We can describe $\mathfrak T_D$ algebraically as follows: To each
element $\mathbf a \in \mathfrak S_{D,n}$ we can associate a
parabolic subgroup of the symmetric group $S_{\mathbf a}$ --- it
is the subgroup preserving the subsets
$\{1,2,\ldots,a_1\},\{a_1+1,\ldots,a_1+a+2\},\ldots,\{D-a_n+1,\ldots,D\}$
of $\{1,2,\ldots,D\}$. Set $T_{\mathbf a}=\sum_{w\in S_{\mathbf
a}}v^{l(w)}T_w$. Then as a module for the Hecke algebra, $\mathfrak T_D$
 is isomorphic to
\[
\bigoplus_{\mathbf a \in \mathfrak S_{D,n}}T_{\mathbf a}\mathcal
H_D
\]
Similarly we see that we can describe an element $[A]$ of $\mathfrak A_D$
uniquely by a triple consisting of an element $w_A \in W$ together with a
pair $\mathbf a,\mathbf b \in \mathfrak S_{D,n}$. Indeed $\mathbf a, \mathbf b$
 are just $r(A)$ and $c(A)$ respectively, and $w_A$ is the element of
maximal length in the (finite) double coset of $S_{\mathbf
a}\backslash W/S_{\mathbf b}$ determined by the matrix $A$. This
also allows us to describe the structure constants for $\mathfrak
A_D$ with respect to the basis $\{[A]\colon A \in \mathfrak S_{D,n,n}\}$
in terms of those for $\mathcal H_D$ with respect to the basis
$\{T_w\colon w \in W\}$. In fact simple algebraic considerations (or an
analogous discussion of the geometry involved) shows that the same
holds for the structure constant with respect to the bases coming
from intersection cohomology.

More precisely, suppose that we denote the various structure
constants for $\mathcal H_D$ and $\mathfrak A_D$ as follows: Let
$A,B \in \mathfrak S_{D,n}$, let $v,w \in W$, and let $\{C_w\colon w
\in W\}$ be the Kazhdan-Lusztig basis of the Hecke algebra.

\begin{itemize}
\item $[A][B] = \sum_C \eta_{A,B}^C[C]$;

\item $\{A\}\{B\} = \sum_C\nu_{A,B}^C\{C\}$;

\item $T_vT_w = \sum_z f_{v,w}^z C_z$;

\item $C_vC_w = \sum_zh_{v,w}^z C_z$.
\end{itemize}
\noindent
Then we have the following relationships between them, which will be crucial in describing
 the cell structure of the affine q-Schur algebra.

\begin{lemma}
Let $A,B,C \in \mathfrak S_{D,n}$, and let $w_A,w_B,w_C \in W$ be
the corresponding element of the Weyl group. Suppose that
$c(A)=r(B)= \mathbf c$. Let $w_{\mathbf c}$ be the longest
element of $S_{\mathbf c}$ and let
$$ p_{\mathbf c}=v^{-l(w_{\mathbf c})}\sum_{x \in S_{\mathbf
c}}v^{2l(x)}$$

be the shifted Poincar\'{e} polynomial of $S_{\mathbf c}$. We have
$$ p_{\mathbf c}\eta_{A,B}^C = f_{w_A,w_B}^{w_C}.$$
$$ p_{\mathbf c}\nu_{A,B}^C = h_{w_A,w_B}^{w_C};$$
\label{coef}
\end{lemma}

\begin{proof}
Using the algebraic description of $\mathfrak T_D$, the first
statement can be proved algebraically, by interpreting the basis
elements $[A],[B],[C]$ as sums of elements in double cosets of the
Hecke algebra. Similarly one can show the second statement
entirely algebraically, but it is perhaps more enlightening to use
the interpretation of the multiplication in terms of perverse
sheaves (see \cite{L99} for a discussion of this). The affine
Schubert variety for $w_A$ fibres over the varieties $\bar
X^{\mathbf L}_A$ with fibre given by a partial flag variety
corresponding to $c(A)$, and the polynomial $p_{\mathbf c}$ arises
from the cohomology of this fibre.
\end{proof}

\bigskip
Finally must describe the connection between our convolution algebras and the modified quantum group
 of affine type A (in its degenerate, or level zero form). We start with some general definitions.
\begin{definition}
A \textit{Cartan datum} is a pair $(I,\cdot)$ consisting of a finite set
$I$ and a $\mathbb Z$-valued symmetric bilinear pairing on the
free Abelian group $\mathbb Z[I]$, such that
\begin{itemize}
\item $i\cdot i \in \{2,4,6,\ldots\}$
\item $2\frac{i\cdot j}{i\cdot i} \in \{0,-1,-2,\ldots\}$, for $i\neq j$.
\end{itemize}
A \textit{root datum} of type $(I,\cdot)$ is a pair $Y,X$ of finitely-generated
 free Abelian groups and a perfect pairing
$\langle,\rangle \colon Y \times X \to \mathbb Z$, together with
imbeddings $I\subset X$, ($i\mapsto i$) and $I\subset Y$, ($i
\mapsto i'$) such that $\langle i,j' \rangle = 2\frac{i\cdot j}{i
\cdot i}$.
\end{definition}

Given a root datum, we may define an associated quantum group $\mathbf U$.
Since it is the only case we need, we will assume that our datum is symmetric
and simply laced so that $i\cdot i = 2$ for each $i \in I$, and $i\cdot j \in \{0,-1\}$
if $i \neq j$. In this case, $\mathbf U$ is generated as an algebra over $\mathbb Q(v)$
by symbols $E_i, F_i, K_\mu$, $i \in I$, $\mu \in Y$, subject to the following relations.
\begin{itemize}
\item $K_0=1$, $K_{\mu_1}K_{\mu_2} = K_{\mu_1+\mu_2}$ for $\mu_1,\mu_2 \in Y$;
\item $K_{\mu} E_i K_{\mu}^{-1} = v^{\langle\mu,i'\rangle}E_i, \quad K_{\mu} F_i K_{\mu}^{-1} =
 v^{-\langle\mu,i'\rangle}F_i$ for all $i \in I$, $\mu \in Y$;
\item $E_iF_j - F_jE_i = \delta_{i,j}\frac{K_i-K_i^{-1}}{v-v^{-1}}$;
\item $E_iE_j=E_jE_i, \quad F_iF_j=F_jF_i$, for $i,j \in I$ with $i\cdot j =0$;
\item $E_i^2E_j +(v+v^{-1})E_iE_jE_i + E_jE_i^2 =0$ for $i,j \in I$ with $i\cdot j = -1$;
\item $F_i^2F_j +(v+v^{-1})F_iF_jF_i + F_jF_i^2 =0$ for $i,j \in I$ with $i\cdot j = -1$.
\end{itemize}

The other object we need is the modified quantum group $\dot{\mathbf U}$.
Let $\mathbf{Mod}_X$ denote the category of left $\mathbf U$-modules
$V$ with a weight decomposition, that is
$$ V = \bigoplus_{\lambda \in X} V_\lambda,$$
where
$$ V_\lambda = \{v \in V \colon K_\mu v = v^{\langle \mu,\lambda\rangle}v, \forall \mu \in Y\}.$$
The forgetful functor to the category of vector spaces has an endomorphism ring $R$. Thus an
element of $a$ of $R$ associates to each $V \in \mathrm{Ob}(\mathbf{Mod}_X)$ an endomorphism $a_V$,
such that for any morphism $f\colon V \to W$, $a_W\circ f = f \circ a_V$.
Thus any element of $\mathbf U$ clearly determines an element of $R$. For each $\lambda \in X$,
let $1_\lambda \in R$ be the projection to the $\lambda$ weight space. Then $R$ is isomorphic to the
direct product $\prod_{\lambda \in X}\mathbf U 1_\lambda$, and we set
$$\dot{\mathbf U} = \bigoplus_{\lambda \in X} \mathbf U 1_\lambda.$$

To see the connection between our convolution algebra and quantum groups,
we will need the following notation. For $\mathbf a \in \mathfrak
S_{D,n}$ let $\mathbf{i_a} \in \mathfrak S_{D,n,n}$ be the
diagonal matrix with $(\mathbf{i_a})_{i,j}=\delta_{i,j}a_i$. Let
$E^{i,j}\in \mathfrak S_{1,n,n}$ be the matrix with
$(E^{i,j})_{k,l}=1$ if $k=i+sn$, $l=j+sn$, some $s \in \mathbb Z$,
and $0$ otherwise. Let $\mathfrak S^n$ be the set of all $\mathbf
b=(b_i)_{i\in \mathbf Z}$ such that $b_i=b_{i+n}$ for all $i \in
\mathbb Z$. Let $\mathfrak S^{n,n}$ denote the set of all matrices
$A= (a_{i,j})$, $i,j \in \mathbb Z$, with entries in $\mathbb Z$
such that
\begin{itemize}
\item $a_{i,j} \geq 0$ for all $i \neq j$; \item
$a_{i,j}=a_{i+n,j+n}$, for all $i,j \in \mathbb Z$; \item For any
$i \in \mathbb Z$ the set $\{j \in \mathbb Z\colon a_{i,j}\neq 0\}$ is
finite; \item For any $j \in \mathbb Z$ the set $\{i \in \mathbb
Z\colon a_{i,j}\neq 0\}$ is finite.
\end{itemize}
Thus we have $\mathfrak S_{D,n,n} \subset \mathfrak S^{n,n}$ for
all $D$. For $i \in \mathbb Z/n\mathbb Z$ let $\mathbf i \in
\mathfrak S^n$ be given by $\mathbf i_k=1$ if $k=i$ mod $n$,
$\mathbf i_k=-1$ if $k=i+1$ mod $n$, and $\mathbf i_k=0$
otherwise. We write $\mathbf a \cup_i \mathbf{a'}$ if $\mathbf
a=\mathbf{a'}+\mathbf i$. For such $\mathbf a,\mathbf{a'}$ set
$_\mathbf a\mathbf e_\mathbf{a'} \in \mathfrak S^{n,n}$ to be
$\mathbf{i_a}-E^{i,i}+E^{i,i+1}$, and $_\mathbf{a'}\mathbf
f_\mathbf a \in \mathfrak S^{n,n}$ to be
$\mathbf{i_{a'}}-E^{i+1,i+1}+E^{i+1,i}$. Note if $\mathbf a,
\mathbf{a'} \in \mathfrak S_{D,n}$ then $_\mathbf a\mathbf
e_\mathbf{a'}, _\mathbf{a'}\mathbf f_\mathbf a \in \mathfrak
S_{D,n,n}$. For $i \in \mathbb Z/n\mathbb Z$ set
\[
E_i(D)= \sum[_\mathbf a\mathbf e_{\mathbf{a'}}], \qquad  F_i(D)=
\sum[_\mathbf{a'}\mathbf f_{\mathbf a}],
\]
where the sum is taken over all $\mathbf a,\mathbf{a'}$ in
$\mathfrak S_{D,n}$ such that $\mathbf a\cup_i\mathbf{a'}$. For
$\mathbf a \in \mathfrak S^n$ set
\[
K_\mathbf a(D) =\sum_{\mathbf b \in \mathfrak S_{D,n}}v^{\mathbf
a\cdot \mathbf b}[\mathbf{i_b}]
\]
where, for any $\mathbf a, \mathbf b \in \mathfrak S^n$, $\mathbf
a\cdot \mathbf b=\sum_{i=1}^n a_ib_i \in \mathbb Z$.
If we let $X'=Y'=\mathfrak S^n$, and $I= \mathbb Z/ n\mathbb Z$, with the embedding
of $I \subset X'=Y'$ and pairing as given above, we obtain a symmetric simply-laced
root datum. We call the quantum group associated to it $\mathbf U(\widehat{\mathfrak{gl}}_n)$.
It can be shown \cite{L99} that the elements $E_i(D), F_i(D), K_{\mathbf a}(D)$,
generate a subalgebra $\mathbf U_D$ which is a quotient of the
quantum group $\mathbf U(\widehat{\mathfrak{gl}}_n)$, via the map the notation suggests.
Note that this gives the algebra $\mathfrak A_D$ the structure of a $\mathbf
U(\widehat{\mathfrak{gl}}_n)$-module. In this paper we will consider the slightly smaller
algebra $\mathbf U = \mathbf U(\widehat{\mathfrak{sl}}_n)$, for which $X = \{\mathbf a \in \mathfrak S^n \colon
 \sum_{i=1}^n a_i = 0 \}$, and $Y = \mathfrak S^n / \mathbb Z \mathbf b_0$, or more precisely
 its modified form which we will denote by $\U$. $\mathbf U(\widehat{\mathfrak{sl}}_n)$ is a subalgebra
 of $\mathbf U(\widehat{\mathfrak{gl}}_n)$, but it is easy to see that its image in $\mathfrak A_D$ is all
 of $\mathbf U_D$.

\section{Distinguished elements in $\mathfrak A_D$}
Our first step in understanding the theory of cells in $\U$ is to
understand the corresponding theory for the affine q-Schur algebra.
We do this by transferring the information known about the Hecke algebra $\mathcal H_D$
 to our case. The key ingredient in our approach is the use of Lusztig's notion of distinguished
elements.

We begin by defining a somewhat mysterious integer-valued function $a_D'$, which
together with certain variants, play a crucial r\^{o}le in our study
of cells.

For $\{C\} \in \mathfrak B_D$ we set $n_{A,B}^C$ to be the largest power of $v$
occurring in the structure constant $\nu_{A,B}^C$, and for
$\mathbf a \in \mathfrak S_{D,n}$ set $|\mathbf a|^2= \sum_{i=1}^n
a_i^2$.

\begin{definition}
For $\ca \in \mathfrak B_D$, such that $\ca = \ca[\mathbf
i_{\mathbf a}]$, consider the set of positive integers
$$\left\{n_{B,C}^A + |\mathbf b|^2 - |\mathbf a|^2\colon \cb, \cc \in \mathfrak
B_D; \cc \in [\mathbf i_{\mathbf b}]\mathfrak A_D[\mathbf i_{\mathbf a}]\right\}.$$
If it has a largest element $d$ we set $a_D'=d$, otherwise
we set $a_D'=\infty$.
\end{definition}

At first sight it would seem that we elided by saying that $a_D'$ is
``integer-valued'' above, however the following lemma shows this is
not the case.
\begin{lemma}
The function $a_D'$ is finite for every $\ca \in \mathfrak B_D$. \label{finite}
\end{lemma}
\begin{proof}
To show this we use the fact that we can interpret the
structure constants in terms of those for the affine Hecke
algebra, and then use the result of Lusztig \cite{L85}, which
shows that the corresponding function on the Kazhdan-Lusztig basis
is finite. Indeed, using Lemma \ref{coef} we see that for $\{A\}, \{B\}, \{C\} \in \mathfrak B_D$
we have $\nu_{A,B}^C = p_{c(A)}^{-1}h_{x,y}^z$ where $x,y,z \in W$
are the corresponding element of the affine Weyl group of type
$\hat{A}_{D-1}$. Now by Theorem 7.2 in \cite{L85} we have
$v^{-l(w_0)}h_{x,y}^z \in \mathbb Z[v^{-1}]$ for any $x,y,z \in
W$, where $w_0$ is the longest element in $S_D$, the finite Weyl
group. The result follows.
\end{proof}

\noindent
Note that we have shown that $a_D'$ is not only finite, but
in fact bounded. We also set $\gamma_{A,B}^C$ to be the coefficient of
$v^{a_D'(C)-|\mathbf b|^2 + |\mathbf a|^2}$ in $\nu_{A,B}^C$
(which in general may be zero).

\begin{definition}
Let $\mathbf a \in \mathfrak S_{D,n}$. For $\ca \in [\mathbf
i_{\mathbf a}]\mathfrak A_D[\mathbf i_{\mathbf a}] $ set
$\Delta(A)$ to be the integer $d \geq 0$ such that
$$ ([\mathbf i_{\mathbf a}],\ca)_D = a_{d}v^{-d} + a_{d+1}v^{-d-1} + \ldots ,$$
where $a_{d} \neq 0$. Set $n_A=a_d$.
\end{definition}

\begin{lemma}
Let $\mathbf a \in \mathfrak S_{D,n}$. For any $\ca \in \mathfrak
B_D$ with $[\mathbf i_{\mathbf a}]\ca[\mathbf i_{\mathbf a}]=\ca$
we have $a_D'(A) \leq \Delta(A)$. \label{ineq}
\end{lemma}
\begin{proof}

This follows an idea of Springer in the Hecke algebra case.
Suppose that $\cb,\cc$ are in $\mathfrak B_D$, and consider the
product $\cb\cc$. We may write this as a sum $\sum_{\{E\} \in
\mathfrak B_D}\nu_{B,C}^E\{E\}$. Chose $\cb \in [\mathbf
i_{\mathbf a}]\mathfrak A_D[\mathbf i_{\mathbf b}]$, and $\cc \in
[\mathbf i_{\mathbf b}]\mathfrak A_D[\mathbf i_{\mathbf a}]$ so
that
$$\nu_{B,C}^{A} = \gamma_{B,C}^{A}v^{a_D'(A)-|\mathbf b|^2 +
|\mathbf a|^2} + \ldots,$$ where $\gamma_{B,C}^{A} \neq 0$ and the
remaining terms are of lower degree. We have the inner product
\begin{equation}
(\cb\cc,[\mathbf i_{\mathbf a}])_D = \sum_{\{E\} \in \mathfrak
B_D}\nu_{B,C}^E (\{E\},[\mathbf i_{\mathbf a}])_D. \label{inp}
\end{equation}

All the terms here have nonnegative integer coefficients (using the positivity of the inner
product). The properties of the  inner product $(,)_D$ show that
this is also equal to
$$v^{|\mathbf a|^2 - |\mathbf b|^2}(\cc,\{B^t\})_D$$
Now since the canonical basis $\mathfrak B_D$ is almost orthogonal with respect to
the inner product, we see that all the terms on the left-hand side
of Equation \ref{inp} lie in $v^{|\mathbf a|^2 - |\mathbf
b|^2}\mathbb N[v^{-1}]$. In particular, taking $\{E\}=\ca$ we
get that
$$ \nu_{B,C}^A(\ca,[\mathbf i_{\mathbf a}])_D=n_A\gamma_{B,C}^A v^{a_D'(A)-|\mathbf b|^2 +
|\mathbf a|^2-\Delta(\ca)}+\ldots  \in v^{|\mathbf a|^2 -|\mathbf
b|^2}\mathbb N[[v^{-1}]],$$

and hence the result.
\end{proof}

Motivated by this, we define the set of \textit{distinguished elements} of
$\mathfrak B_D$ as follows.

\begin{definition}
Let $\mathcal D_D$ be the set of elements $\ca$ in $\mathfrak B_D$
such that there is a $\mathbf a \in \mathfrak S_{D,n}$ with $\ca
\in [\mathbf i_{\mathbf a}]\mathfrak A_D[\mathbf i_{\mathbf a}]$
and $a_D'(A)=\Delta(A)$.
\end{definition}

The distinguished elements $\mathcal D_D$ are defined by analogy with
 the Hecke algebra case due to Lusztig \cite{L87}. We note the some
 consequences of the above proof.

\begin{cor}
We have the following properties:
\begin{enumerate}
\item If $\ca \in \mathcal D_D$ and $\cb,\cc \in \mathfrak B_D$
are such that $\gamma_{B,C}^A \neq 0$ then $\cc=\{B^t\}$. \item
For each $\cb \in \mathfrak B_D$, there is a unique $\ca \in
\mathcal D_D$ with $\gamma_{B,B^t}^A \neq 0$ \item If $\ca \in
\mathcal D_D$ then $\ca=\{A^t\}$.
\end{enumerate}
\label{list}
\end{cor}
\begin{proof}
For the first, note that in the above proof, the almost
orthogonality of $\mathfrak B_D$ with respect to the inner product
implies that it is necessary and sufficient to have $\cc =
\{B^t\}$. That the product contains just one element of $\mathcal
D_D$ is also immediate. For the last statement, pick $\cb,\cc$
such that $\gamma_{B,C}^A \neq 0$. By the first statement, we see
that $\cc=\{B^t\}$. Since the product $\cb\{B^t\}$ is preserved by
the transpose anti-automorphism $\Psi$, we see that $\gamma_{B,B^t}^{A^t}
\neq 0$, and so by the second statement, $\ca = \{A^t\}$.
\end{proof}

Recall that to each element of $\mathfrak B_D$ we have attached an
element of the affine Weyl group. We will show that in this way,
the distinguished elements of $\mathfrak  B_D$ actually correspond
to distinguished elements of W.

\begin{lemma}
Let $\ca \in \mathcal D_D$, then the Weyl group element $w_A$ is
distinguished and conversely.
\label{Duflo}
\end{lemma}
\begin{proof}
Let $\mathbf a \in \mathfrak S_{D,n} $ be such that $[\mathbf
i_{\mathbf a}]\ca = \ca$. By definition we see that
$$(\ca,[\mathbf i_{\mathbf a}])_D= \Pi_{\mathbf i_{\mathbf a},E},$$
the stalk of the intersection cohomology sheaf on
$\bar{X}^{\mathbf L}_E$ at the point corresponding to
$[\mathbf{i}_{\mathbf a}]$. But this is equal to $v^{l(w_{\mathbf
a})}p_{1,w_E}$, where $p_{1,w_E}$ is the affine Kazhdan-Lusztig
polynomial attached to $1,w_E$, and $w_{\mathbf a}$ is the longest
element of the parabolic subgroup attached to $\mathbf a$ (the intersection
cohomology sheaves are related by a smooth pullback with fibre
dimension $l(w_{\mathbf a})$. Thus we see that if $\Delta(w_E)$ is the lowest
power of $v^{-1}$ occurring in $p_{1,w_E}$,
$$a_D'(E) \leq a'(w_E)-l(w_{\mathbf a}) \leq \Delta(w_E)
-l(w_{\mathbf a})=\Delta(E).$$ Here the function $a'$ on the
Kazhdan-Lusztig basis is the one defined in \cite{L87} (there
denoted simply $a$). For $z \in W$, we set $a'(z)$ to be the highest
power of $v$ appearing in a structure constant $h_{x,y}^{z}$ as
$x,y$ vary over $W$. The first inequality follows directly from
the definitions of $a'$, $a_D'$, and the second from analog of Lemma
\ref{ineq} for $\mathcal H_D$. It follows immediately that $w_E$
is distinguished if $\{A\}$ is. To establish the converse, it is necessary
to note that if one picks any element $x$ of the left cell
containing the distinguished element, then by \cite{L87} the structure constant
$h_{x^{-1},x}^{w_E}$ has $a'(w_E)$ as its highest power of
$v$, and so $a_D'(E)=a'(w_E)-l(w_{\mathbf a})$.
\end{proof}

We now show that all the notions of cells for $\mathfrak A_D$ can
be deduced from those for the Hecke algebra. More precisely we
have the following result.

\begin{prop}
Let $\ca,\cb \in \mathfrak B_D$.
\begin{enumerate}
\item $\ca \sim_L \cb$ if and only if $w_A \sim_L w_B$ and
$c(A)=c(B)$;
\item $\ca \sim_R \cb$ if and only if $w_A \sim_R w_B$ and $r(A)=r(B)$;
\item $\ca \sim_{LR} \cb$ if and only if $w_A \sim_{LR} w_B$;
\item $a_D'(\ca)=a'(w_A)-l\left(w_{c(A)}\right)$;
\item Each left cell contains precisely one distinguished element.
\end{enumerate}
\end{prop}
\label{cellprop}
\begin{proof}
For the first claim, note that since the notion of cell in
$\mathfrak A_D$ is defined essentially by using a subset of the
Kazhdan-Lusztig basis consisting of those elements which are of
maximal length in certain double cosets, it is clear that if $\ca
\sim_L \cb$ then $w_A\sim_Lw_B$. Moreover, certainly we have
$c(A)=c(B)$. For the converse, we need to use the distinguished
elements. Suppose that $w_A\sim_L w_B$ and $c(A)=c(B)$. Then if
$d$ is the unique distinguished element in the left cell $\Gamma$
containing $w_A,w_B$ (which exists by \cite{L87}), $d$ determines
a distinguished element of $\mathfrak B_D$, $\{E\}$ say, where
$c(E)=c(A)$.

For $x,y,z \in W$ let $\gamma_{x,y}^z$ denote the coefficient of $v^{a'(z)}$
in $h_{x,y}^z$. Then by \cite[Theorem 1.8]{L87} we know that
$\gamma_{x,y}^z=\gamma_{y,z^{-1}}^{x^{-1}}=\gamma_{z^{-1},x}^{y^{-1}}$, and
so as $\gamma_{w_A^{-1},w_A}^{d} \neq 0$ we see that
$$h_{w_A^{-1},w_A}^{d},h_{w_A^{-1},d}^{w_A},h_{w_B^{-1},w_B}^{d},h_{w_B^{-1},d}^{w_B}$$
are all nonzero, and hence the same is true of
$$\nu_{A^t,A}^E,\nu_{A^t,E}^A,\nu_{B^t,B}^E,\nu_{B^t,E}^B.$$ It
follows that $\ca \sim_L \{E\}$ and $\cb \sim_L \{E\}$, and hence
$\ca \sim_L \cb$.

The second claim either follows in the same way, or by taking inverses in $W$,
which corresponds to applying the transpose map $\Psi$ in $\mathfrak A_D$.

For the third, the forward implication is again clear. If $w_A
\sim_{LR} w_B$, then it follows that the left cell containing
$w_A$ and the right cell containing $w_B$ intersect (since this is
true of any left and right cell in the same two-sided cell of an
affine Hecke algebra). As any element in this intersection will
give rise to an element $\cc$ of the q-Schur algebra with
$r(C)=r(B)$ and $c(C)=c(A)$, we obtain the result using the first
two parts of the proposition. The fourth claim follows from Corollary
\ref{list} (see the end of the proof of Lemma \ref{Duflo} . Since
each left cell of the Hecke algebra contains a unique distinguished element,
the fifth claim follows from the Lemma \ref{Duflo} and the first claim.
\end{proof}

\section{Cells in $\mathfrak A_D$}

We saw at the end of the last section that the theory of cells of the
affine q-Schur algebra is determined by that for the type A affine Hecke
algebra. This allows us to describe explicitly the number of two-sided cells in
the affine q-Schur algebra, and also the number of left cells (and
hence right cells) in a given two-sided cell. To do this we
recall the combinatorial definitions which describe the bijection
between cells for the Hecke algebra and partitions.

\begin{definition}
Suppose $w \in W$ the affine Weyl group of type $\tilde A_{D-1}$.
Then we may view $w$ as a permutation of $\mathbb Z$. A sequence
$(i_1,i_2,\ldots,i_r)$ is called a \textit{d-chain} if $i_1<i_2< \ldots
<i_r$ and $(i_1)w>(i_2)w>\ldots>(i_r)w$, and the $\{i_j,
j=1,\ldots,r\}$ are all incongruent modulo $D$.

Let $\mathcal P_D$ be the set of partitions of $D$.
We define a map $\sigma\colon W \to \mathcal P_D$ as follows. For $w \in W$,
let $d_j$ be the maximal size of a set of $j$ d-chains, the union of whose
elements are all incongruent modulo $D$. Then it is know that
$\lambda=(d_1,d_2-d_1,d_3-d_2,\ldots,d_{D-1}-d_D)$ is a partition
of $D$. Set $\sigma(w)= \lambda$.
\end{definition}

The following result is due to Lusztig, based on the work of Shi, see \cite{L85a}, \cite{Sh}.
\begin{theorem}
The fibres of $\sigma$ are precisely the two-sided cells of $W$.
\hfill $\square$ \label{part}
\end{theorem}

We may use this to give a description of the two-sided cells in
the affine q-Schur algebra as follows:

\begin{definition}
Let $A \in \mathfrak S_{D,n}$. An \textit{anti-diagonal path} in $A$ is an infinite
strip of entries $(a_{i_k,j_k}\colon k \in \mathbb Z)$ such that $(i_k,j_k)$ is either
equal to $(i_{k-1}-1,j_{k-1})$ or $(i_{k-1},j_{k-1}+1)$ with the latter being the
case for all but finitely many $k$. Thus visually if you draw the matrix with rows
increasing from top to bottom, and columns from left to right, (as we will do)
then path starts and ends with infinite vertical strips, and takes finitely many right or
vertical turns.

Let $d_j$ be the maximal size of the sum of entries in the union of $j$ anti-diagonal paths.
Then we define a map $\rho\colon \mathfrak S_{D,n} \to \mathcal P^n_D$ where
$\mathcal P^n_D$ is the set of partitions of $D$ with at most n parts,
by setting $\rho(A)=(d_1,d_2-d_1,\ldots,d_n-d_{n-1})$. As above, it
follows from general results on posets that $\rho(A)$ is indeed a partition.
The fact that it can have at most $n$ parts is obvious. We will sometimes
view $\rho$ as a map from $\mathfrak B_D$ in the obvious way.
\end{definition}

\begin{prop}
The fibres of the map $\rho\colon \mathfrak B_D \to \mathcal P^n_{D}$ are the two-sided cells
of $\mathfrak A_D$.
\label{partqSchur}
\end{prop}
\begin{proof}
This is a simple combination of the statements of Proposition
\ref{cellprop} and Theorem \ref{part}.
\end{proof}

\begin{example}
Suppose that $n=2$ and $D=5$. Consider the element $\ca$ of $\mathfrak B_5$ corresponding to
\[
\left(
\begin{array}{cccccccc}

& \vdots & \vdots &\vdots & \vdots & \vdots & &  \\
\cdots& 1& 1& 0& \fbox{1}& 0& \cdots \\
\cdots& 0& 1& 0& \fbox{1}& 0& \cdots \\
\cdots& 0& 0& \fbox{1}& \fbox{1}& 0& \cdots \\
\cdots& 0& 0& \fbox{0}& 1& 0& \cdots \\
& \vdots & \vdots & \vdots & \vdots & \vdots & &  \\

\end{array}
\right)
\]
where the top left entry shown is in the $(1,1)$ entry of $A$. Then $\ca$ lies in the two-sided cell
corresponding to the partition $(4,1)$. The boxed entries give part of an anti-diagonal path
which has entry sum $4$. Note that it is not unique.
\end{example}

Given a partition $\lambda \in \mathcal P^n_D$ we denote the two-sided cell $\rho^{-1}(\lambda)$ by
$\mathbf c_\lambda$. We will often use the same notation for a partition in $\mathcal P_D$ and a two-sided
cell of $\mathcal H_D$.

Somewhat more elaborate is a description of the number of left
cells in a two-sided cell of the affine q-Schur algebra. Notice that Proposition
\ref{cellprop} shows that each left cell of the affine Hecke
algebra gives rise to a number of left cells of the affine q-Schur
algebra, with the number depending on the set of simple
reflections of the symmetric group $S_D$ which decrease the length
of any element of the left cell when multiplied on the right.

For $w \in W$ let $\mathcal R (w)= \{s \in S\colon l(ws) < l(w)\}$ and
$\mathcal L (w)= \{s \in S\colon l(sw) < l(w)\}$.
It is known that the functions $\mathcal R, \mathcal L$ are constant on
right and left cells respectively. Thus for $\Gamma$ a left cell, we may
write $\mathcal R(\Gamma)$ for the set $\mathcal R(w)$, where $w$ is any
element of $\Gamma$. The left cells of $\mathcal H_D$ have been described by
Shi (\cite{Sh}, chapter 14) as the fibres of a map to a set of tableaux, such that
the shape of the tableau associated to a left cell is given by the partition of the
two-sided cell it lies in, and the entries must increase down the columns.

In order to describe this map in more detail we need some to make some definitions.
We use the description of $W$ as a group of permutations of the integers,
and in particular the associated infinite matrices. Let $A= (a_{i,j})_{i,j \in \mathbb Z}$
be such a matrix. A \textit{block} is a set of consecutive rows of $A$.
For a block of $m$ rows $i+1,i+2,\ldots,i+m$, let the nonzero entries be
$\{a_{i+1,j_1},a_{i+2,j_2},\ldots,a_{i+m,j_m}\}$. We say the block is a \textit{descending chain} if
$j_1 > j_2 > \ldots > j_m$. A block is a \textit{maximal descending chain} (MDC) if it cannot be
imbedded in a larger such block.

Say that an element of $w \in W$ has \textit {full MDC form} at $i$, if there
exist consecutive MDC blocks $(A_l,A_{l-1},\ldots,A_1)$ of $A_w$ of size $m_t$, for $t=1,\ldots,l$,
with $\sum_{t=1}^l m_t = D$, and $i+1$ the first row of $A_l$ ( so $i+\sum_{j=1}^k m_j + 1$ is the
first row of $A_{l-k}$).
Suppose a full MDC form has blocks which are of (weakly) increasing size (so $A_i$ has at most as many
rows as $A_{i-1}$). Let $j_t^u$ be the column of the nonzero entry in the $u$-th row of $A_t$.
Then we say the form is \textit{normal} if $j_1^u-n < j_l^u < j_{l-1}^u < \ldots < j_1^u$, for each
$u$ (where we ignore terms in this sequence which do not exist).

For $\lambda=(\lambda_1 \geq \lambda_2 \geq \ldots \geq \lambda_r >0)$ a partition of $D$, let
$N_\lambda$ be the set of elements of the two-sided cell $\mathbf c_{\lambda}$ corresponding
to $\lambda$ which have normal MDC form $(A_r,A_{r-1},\ldots,A_1)$ at $i$ for some $i \in \mathbb Z$,
where $\lambda_j$ is the number of rows in $A_j$.

\begin{theorem}\cite{Sh}
Let $w \in \mathbf c_\lambda$. Then there is a $y \in N_\lambda$ with $y \sim_L w$.
\hfill $\square$
\end{theorem}

Let $\mathcal C_\lambda$ be the set of Young diagrams of shape $\lambda$ with entries $\{1,2,\ldots,D\}$
which decrease down columns. In \cite[chapter 14]{Sh} Shi defines a map $T$ from the left cells
in $\mathbf c_\lambda$ to $\mathcal C_\lambda$. Let $\Gamma$ be a left cell in $\mathbf c_\lambda$.
Choose $y \in N_\lambda \cap \Gamma$ and then set the entries of column $u$ of $T(\Gamma)$ to be the residues modulo
$D$ of the numbers $\{j_t^u\colon 1 \leq t \leq \mu_u\}$ where $\mu$ is the partition dual to $\lambda$.
Shi shows this is independent of the choice of $y$, and that it gives a bijection.

\begin{example}
Consider the matrix $A_w$ in $\hat{A}_4$ given as follows
\[
\left(
\begin{array}{ccccccc}
\cdots & 0 & 0 & 1 & 0 & 0 & \cdots \\
\cdots & 0 & 1 & 0 & 0 & 0 & \cdots \\
\cdots & 0 & 0 & 0 & 0 & 1 & \cdots \\
\cdots & 0 & 0 & 0 & 1 & 0 & \cdots \\
\cdots & 1 & 0 & 0 & 0 & 0 & \cdots \\
\end{array}
\right)
\]
where the first column  shown is column 1. Then it is easy check that $w \in N_{(3,2)}$ and the associated
tableau is given below.

\[
\begin{tableau}
:.5.4.1 \\
:.3.2 \\
\end{tableau}
\]
\end{example}

It follows directly from this construction (though this is not explicitly
described in \cite{Sh}) that the set of simple reflections in $\mathcal R(\Gamma)$
is determined by this tableau. Indeed since $\mathcal R(\Gamma)$ is given
by $\mathcal R(w)$ for any $w \in \Gamma$, we may use an element of $N_\lambda$
as above. Then it is easy to see that the simple reflection $s_i$ is in
$\mathcal R(\Gamma)$ precisely when $i$ appears to the right of $i+1$ in the
tableau (where one reads modulo $D$ for $s_0$).

We now consider the left cells of $\mathfrak A_D$. Each such cell correspond to a left
cell $\Gamma$ of $\mathcal H_D$ and an element $\mathbf a$ of $\mathfrak S_{D,n}$
such that the simple reflections $J$ of $S_{\mathbf a}$ are a subset of
$\mathcal R(\Gamma) \backslash \{s_0\}$.

\begin{definition}
For $\lambda \in \mathcal P_D$, let $\mathcal C^n_{\lambda}$ be the
the set of tableaux of shape $\lambda$ with entries from $\{1,2,\ldots,n\}$ strictly decreasing
down columns. Thus $\mathcal C^n_{\lambda}$ is empty if $\lambda$ has more than $n$ parts.
\end{definition}

If we fix a two-sided cell $\mathbf c_\lambda$, where $\lambda \in \mathcal P^n_D$,
using the description of $\mathcal R(\Gamma)$ in terms of the tableau $T(\Gamma)$, it is easy
to see that the left cells in $\mathbf c_\lambda$ are indexed by the elements of
$\mathcal C^n_{\lambda}$. Indeed to each tableau $T \in \mathcal C^n_{\lambda}$ there is a
well-defined element $h(T)$ of $\mathcal C_\lambda$ given as follows.
Order the boxes of $T$ by listing those labelled $1$ first, then $2$, and so on, always reading
from right to left. Then construct $h(T) \in \mathcal C_\lambda$ by labelling each box with
its position in the order just described. This gives the left cell of $W$. The element $\mathbf a$
is determined by letting $a_i$ be the number of boxes of $T$ labelled $i$, for $i \in \{1,2,\ldots,n\}$.
The following example makes the correspondence clear.

\begin{example}
Let $D=5$ and $n=3$. Suppose that we consider the tableau
\[
\begin{tableau}
:.3.2\\
:.2.1\\
:.1\\
\end{tableau}
\]
in $\mathcal C^3_{(2,2,1)}$. Then the tableau corresponding to it in $\mathcal C_{(2,2,1)}$ is
\[
\begin{tableau}
:.5.3\\
:.4.1\\
:.2\\
\end{tableau}
\]
and the sequence $\mathbf a$ is $(2,2,1)$ (repeated periodically).
\end{example}

This allows us to count the number of left cells in a two-sided cell of $\mathfrak A_D$.
\begin{prop}
Let $\mathbf c$ be a two-sided cell of $\mathfrak A_D$. If
$\lambda$ is the partition of $D$ associated to $\mathbf c$, and
$\lambda(i)\colon=\lambda_i-\lambda_{i+1}$, then the number of left
cells in $\mathbf c$ is
$$\prod_{i=1}^{n-1}
\binom{n}{i}^{\lambda(i)}$$
\hfill $\square$
\label{leftcellnumber}
\end{prop}

Finally we wish to construct the asymptotic algebra associated to a
two-sided cell. We will need a variant of the function $a_D'$.

\begin{definition}
Let $\ca \in \mathfrak B_D$. If there is an integer $d \geq 0$
such that $v^{-d}\nu_{A,B}^{C} \in \mathbb Z[v^{-1}]$ for all
$\cb,\cc \in \mathfrak B_D$ then let $a(A)$ be the smallest such.
Otherwise set $a(A)=\infty$.
\end{definition}

Note that the proof of Lemma \ref{finite} shows that $a_D$ is always finite.
More interestingly we have the following result.

\begin{lemma}
The functions $a_D'$ and $a_D$ agree. Moreover the function $a_D$ is constant
on $\mathbf c[\mathbf i_{\mathbf a}]$ for any two-sided cell $\mathbf c$ and
any $\mathbf a \in \mathfrak S_{D,n}$.
\label{aconstant}
\end{lemma}
\begin{proof}
Both of these follow from facts about the Hecke algebra: If
$\gamma_{x,y}^z$ denotes the coefficient of $v^{a'(z)}$ in
$h_{x,y}^z$, then it follows from the results above that for
$\ca,\cb,\cc \in \mathfrak B_D$ we have
$\gamma_{A,B}^C=\gamma_{w_A,w_B}^{w_C}$. Moreover, by the results
of \cite{L87} (see Lemma \ref{cellprop}) we know that $\gamma_{x,y}^z=\gamma_{y,z^{-1}}^{x^{-1}}=
\gamma_{z^{-1},x}^{y^{-1}}$ and hence
$\gamma_{A,B}^C=\gamma_{B,C^t}^{A^t}=\gamma_{C^t,A}^{B^t}$. It is
now easy to see that $a_D'=a_D$. Since $a'$ is constant on
two-sided cells of the Hecke algebra, the second statement is
clear.
\end{proof}

We now rescale the canonical basis of $\mathfrak A_D$. Set $\langle A \rangle=v^{-a_D(A)}\ca$. Let
$\mathfrak A_{\mathbf c}$ denote the span of the elements in $\mathbf c$ a two-sided cell of $\mathfrak B_D$.
This becomes an algebra by identifying it with a subquotient of $\mathfrak A_D$ in the obvious way.
The structure constants of $\mathfrak A_{\mathbf c}$ with respect to the new basis lie in $\mathbb Z[v^{-1}]$.
Indeed the product $\langle A \rangle \langle B \rangle = \sum_{\langle C\rangle} v^{-a_D(A)} \nu_{A,B}^C \langle C \rangle$,
since $a_D(A)=a_D(C)$, and the coefficients $v^{-a_D(A)} \nu_{A,B}^C$ all lie in $\mathbb Z[v^{-1}]$.
Thus if $\mathcal L_{\mathbf c}$ is the $\mathbb Z[v^{-1}]$ span of the $\{\langle A \rangle\colon \{A\} \in \mathbf c\}$,
$\mathcal L_{\mathbf c}$ has the structure of a $\mathbb Z[v^{-1}]$ algebra. The quotient $J_{\mathbf c}=
\mathcal L_{\mathbf c}/v^{-1}\mathcal L_{\mathbf c}$ is then a $\mathbb Z$
algebra, where if $t_A$ is the image of $\langle A \rangle$, the multiplication in $J_{\mathbf c}$ is given by
$$t_At_B = \sum_{C}\gamma_{A,B}^{C}t_C.$$
Let $\mathcal D_{\mathbf c} = \mathcal D_D\cap \mathbf c $. It follows from the above that the set
 $\{t_E\colon \{E\} \in \mathcal D_{\mathbf c} \}$ gives a decomposition of the identity into orthogonal idempotents.

By using the results of \cite{Xi} or \cite{BO} we can also give an explicit description of this asymptotic
 algebra. For $\lambda \in \mathcal P^n_D$ and $i \in \{1,2,\ldots,n\}$, let $\lambda(i) = \lambda_i - \lambda_{i+1}$,
 (where $\lambda_{n+1}=0$). Let $G_\lambda$ be the reductive group $\prod_{i=1}^n GL_{\lambda(i)}(\mathbf C)$
 and let $R_\lambda$ be the $K$-group of its representations, so that the irreducible representations
 $\widehat G_\lambda$ form a $\mathbb Z$-basis of $R_\lambda$. Let $T_\lambda$ be the set of triples $(E_1,E_2,\kappa)$
  where $\{E_1\},\{E_2\} \in \mathcal D_\lambda$, and $\kappa \in \widehat{G}_\lambda$. Let $\mathcal J_\lambda$
 be the free Abelian group on $T_\lambda$. Define a ring structure on $\mathcal J_\lambda$ by
$$(E_1,E_2,\kappa)(E_1'E_2',\kappa')=\sum c_{\kappa,\kappa'}^{\kappa''}\delta_{E_2,E_1'}(E_1,E_2',\kappa'')$$
where the sum is over $\kappa'' \in \widehat{G}_\lambda$ and $c_{\kappa,\kappa'}^{\kappa''}$ is the multiplicity
of $\kappa''$ in the $G_\lambda$-module $\kappa\otimes\kappa'$. Thus $\mathcal J_\lambda$ is a matrix ring of rank $N$ over
the representation ring $R_\lambda$, where $N$ is the number of left cells in $\mathbf c_\lambda$ given in Proposition
 \ref{leftcellnumber}.

\begin{prop}
\begin{enumerate}
\item There is a ring isomorphism $J_{\mathbf c_\lambda} \to \mathcal J_\lambda$ which restricts to a bijection between the
canonical basis of $J_{\mathbf c_\lambda}$ and $T_\lambda$.
\item For any $\{E\} \in \mathcal D_{\mathbf c_\lambda}$, the subset of $\mathbf c_\lambda$ corresponding to
$\{(E_1,E_2,\kappa) \in T_\lambda\colon E_2=E\}$ under the bijection is a left cell.
\item For any $\{E\} \in \mathcal D_{\mathbf c_\lambda}$, the subset of $\mathbf c_\lambda$ corresponding to
$\{(E_1,E_2,\kappa) \in T_\lambda\colon E_1=E\}$ under the bijection is a right cell.
\end{enumerate}
\label{Jstructure}
\hfill $\square$
\end{prop}

This shows that all the simple modules of the $\mathbb C$-algebra $\mathbb C \otimes J_{\mathbf c_\lambda}$
 are $N$-dimensional and that the set of isomorphism classes of such modules is in bijection with the semisimple
 conjugacy classes of $G_\lambda$.

The asymptotic algebra also receives a homomorphism from the original algebra, once we tensor
with $\mathbb Q(v)$. Define a map $\Phi_{\mathbf c_\lambda}\colon \mathfrak A_D \to \mathbb Q(v)\otimes J_{\mathbf c_\lambda}$
 as follows:
$$ \Phi_{\mathbf c_\lambda}(\{A\}) = \sum_{\{E\} \in \mathcal D_{\mathbf c_\lambda}, \{B\} \in \mathbf c_\lambda}\nu_{A,E}^{B}t_B.$$
One shows it is a homomorphism as in \cite[Proposition 1.9]{L95}, where the property of the structure
constants which is needed follows from the Hecke algebra case. This allows one to pull back representations
of $J_{\mathbf c_\lambda}$ to representations of $\mathfrak A_D$.

\section{Cells in $\U$}
\label{cellsU}

Let $\U$ be the modified quantum group of affine $\mathfrak{sl}_n$, and let $\B$ be its canonical basis
 (see \cite{L93}). We now show how we can lift information about the cell structure of the affine q-Schur
 algebra to the modified quantum group. If $\phi_D$ were surjective this would be a straightforward consequence
 of the previous section. Indeed in the finite type case (see \cite{BLM90}, \cite{L99a}), the analogue of
 $\phi_D$ is surjective, and the results of \cite{L95} in the case of $\mathfrak{sl}_n$ can be recovered in this way,
 as was essentially done by Du in \cite{Du} (note however that \cite{L95} is much more general, classifying the cell
 structure for any finite type quantum group).

In the affine case it is no longer true that the homomorphism from the quantum group is surjective.
 Thus we need to be more careful in lifting information from $\mathfrak A_D$ to $\U$. The following theorem relating
 the canonical bases $\dot{\mathbf B}$ and $\mathfrak B_D$ was conjectured by Lusztig, and proved in \cite{ScV}.
 A more geometric proof can be found in \cite{M02}.

\begin{theorem}
For all $b \in \dot{\mathbf B}$ we have $\phi_D(b)
\in \{0\}\cup\mathfrak B_D$. Moreover the kernel of $\phi_D$ is
spanned by the elements $b \in \dot{\mathbf B}$ such that
$\phi_D(\mathbf b)=0$. \label{compatible}
\hfill $\square$
\end{theorem}

It follows that the image $\mathbf U_D$ is a union of two-sided cells of $\U$. Moreover
the injectivity result of \cite{L99a} (see also \cite{M02}) shows that any two-sided
cell will eventually lie in some $\mathbf U_D$. Given $A \in \mathfrak S_{D,n}$ we say that $A$
is \textit{aperiodic} if, for any integer $k \neq 0$ there is an integer $p$ with $a_{p,p+k}=0$. Thus
only the main diagonal of $A$ can consist entirely of nonzero entries.
In \cite{L99}, Lusztig showed that $\mathbf U_D$ is spanned by a subset $\mathbf B_D$ of $\mathfrak B_D$
consisting of those $\{A\}$ for which $A$ is aperiodic.

We now define an analogue of the $a_D$ function, following \cite{L95}.
Let $c_{b,b'}^{b''}$ be the structure constants of $\U$ with respect
to $\B$. For a two-sided cell $\mathbf c$ in $\U$ let $\U_{\mathbf c}$
be the subspace of $\U$ spanned by the elements of $\mathbf c$. We endow
$\U_{\mathbf c}$ with an algebra structure by identifying it with a
subquotient of $\U$, so that for $b,b' \in \mathbf c$ the product is given by

$$ bb'= \sum_{b'' \in \mathbf c}c_{b,b'}^{b''}b''.$$

\begin{definition}
Let $b \in \B$. If there is an integer $n \geq 0$ such that
$v^{-n}c_{b,b'}^{b''} \in \mathbb Z[v^{-1}]$ for all $b',b'' \in
\mathbf c$ then let $a(b)$ be the smallest such. Otherwise
set $a(b)=\infty$.
\end{definition}

The following observation simple observation tells us about the left cells
 in $\U$.

\begin{lemma}
Let $\Gamma$ be a left cell of $\mathfrak A_D$, and let $\{E\} \in
\mathcal D_D$ be the unique distinguished element in $\Gamma$.
Then if $\mathbf B_D\cap\Gamma \neq \emptyset$ we must have $\{E\}
\in \mathbf B_D$. Moreover $\mathbf B_D\cap\Gamma$ is a single
left cell of $\mathbf U_D$
\label{leftcell}
\end{lemma}
\begin{proof}
Pick $\ca \in \mathbf B_D\cap\Gamma$. Then we know that
$$\{A^t\}\ca = \nu_{A^t,A}^E\{E\} + \dots,$$
where $\nu_{A^t,A}^E \neq 0$ since $\gamma_{A^t,A}^E \neq 0$ (the
unique distinguished element for which $\gamma_{A^t,A}^E \neq 0$ must
clearly be the one in the left cell containing $\ca$). This
implies that $\{E\} \in \mathbf B_D$. By arguing as in the proof
of the first claim in Proposition \ref{cellprop}) we see that the
intersection $\mathbf B_D\cap\Gamma$ is a single left cell.
\end{proof}

This has some important corollaries which we now record.

\begin{cor}
We have the following properties of $a$ functions.
\begin{enumerate}
\item The functions $a_D,a_D'$ coincide with the analogous functions defined in terms of
$\mathbf U_D$ instead of $\mathfrak A_D$.
\item For $b \in \B$ if $\phi_D(b) \neq 0$ then $a(b)=a_D(\phi_D(b))$. In particular, $a(b)$ is finite.
\end{enumerate}
\label{finitecells}
\end{cor}
\begin{proof}
The claims are easy consequences of the above lemma, using distinguished elements.
\end{proof}

We now know that each left cell in $\U_D$ contains a unique distinguished element.
We wish to show that the notion of distinguished elements lifts to $\U$. Since any
left cell will occur as a left cell of $\U_D$ for sufficiently large $D$, it suffices
to show that the distinguished element we obtain is independent of $D$. Since the distinguished
element is characterized as the idempotent element in the asymptotic algebra, which is determined
by the two-sided cell, it is independent of the algebra $U_D$ we choose. Moreover one of the main
results of \cite{M02} is that the inner product on $\U$ is obtained as a limit from those on $\mathfrak A_D$,
 thus we may give an intrinsic characterization of the set $\mathcal D_{\mathbf c}$ of distinguished elements in a
two-sided cell $\mathbf c$. Recall from \cite[3.7]{L95} that $\U$ possesses an anti-automorphism
$\sharp \colon \U \to \U$, which is such that $\phi_D(x^\sharp) = \Psi(\phi_D(x))$, for any $x \in \U$.

\begin{prop}
Let $b \in \mathbf c$ and $\lambda \in X$ be such that $b \in \U 1_{\lambda}$. Then
$v^{a(b)}(1_{\lambda},b) \in \mathbb Z[v^{-1}]$ with nonzero constant term for
$b \in \mathcal D_{\mathbf c}$ and $v^{a(b)}(1_{\lambda},b) \in v^{-1}\mathbb Z[v^{-1}]$
otherwise. Moreover $b = b^{\sharp}$. \hfill $\square$
\end{prop}

It remains to investigate the structure of the two-sided cells of $\mathbf U_D$. This again lifts from $\mathfrak A_D$,
 as the following simple observations show: The transfer map $\psi_D\colon \mathbf U_D \to \mathbf U_{D-n}$ is such that
 $\psi_D(\ca) = \{A-I\}$ if the entries of $A-I$ are nonnegative and $\psi_D(\ca) =0$ otherwise ($I = (\delta_{ij})$
 is the identity matrix). Using this along with our combinatorial description of two-sided cells in $\mathfrak A_D$,
 we obtain the following statement. Let $\lambda = (\lambda_1 \geq \lambda_2 \geq \ldots \geq \lambda_n \geq 0)$
 be in $\mathcal P^n_D$, and let $\mathbf k_\lambda$ denote the intersection $\mathbf B_D \cap \mathbf c_\lambda$.
 Then $\mathbf k_\lambda$ is a union of two-sided cells of $\mathbf U_D$ and moreover it follows from the above
 discussion that $\mathbf k_\lambda$ maps to $0$ under $\psi_D$ unless $\lambda_n >0$, when it maps to
 $\mathbf k_{\lambda'}$ in $\mathbf U_{D-n}$ where $\lambda'= (\lambda_1-1,\lambda_2-1,\ldots,\lambda_n-1)$
 (i.e.$\lambda'$ is obtained by removing the first column of the Young diagram for $\lambda$).

However, the following observation which follows easily from Proposition \ref{partqSchur} now shows that
we are almost done.

\begin{lemma}
Let $A$ be an element of $\mathfrak S_{D,n}$. Then $\rho(A)$ has strictly less than $n$ parts, precisely when
$A$ has no completely nonzero diagonal (i.e. for each $k \in \mathbb Z$ there is some $p \in \mathbb Z$ with
$a_{p,p+k}=0$). In particular, if $\lambda \in \mathcal P^n_D$ has fewer than $n$ parts
$\mathbf c_\lambda$ consists entirely of aperiodic elements.
\hfill $\square$
\end{lemma}
\noindent
Thus for such $\lambda$ we see that $\mathbf k_\lambda= \mathbf c_\lambda$, and it consists of a single two-sided
cell of $\mathbf U_D$, or $\U$.

Recall the group $X$ of the root datum of $\U$ from section \ref{qschur}. For convenience, here we will view
it as a quotient of $\mathbb Z^n$ (by taking the entries $a_1,a_2,\ldots,a_n$).
We define $X^+$ to be the``dominant weights'' in $X$. Let $I_0 = \{i \in
\mathbb Z/n \mathbb Z\colon i \neq 0 $ mod $ n \}$. The set $X^+$ consists of those $\mu \in X$
with $\mu(i) \colon=\langle i,\mu \rangle \geq 0$ for $i \in I_0$.

\begin{prop}
The two-sided cells of $\U$ are naturally parameterized by $X^+$.
\end{prop}
\begin{proof}
First note that each partition $\lambda$ with at most $n$ parts determines an element $\overline{\lambda}$ in $X^+$ by
 taking the coset of $(\lambda_1,\lambda_2,\ldots,\lambda_n)$ in $X$, and the previous paragraph
 shows that this gives a natural bijection between $X^+$ and the two-sided cells of $\U$. Indeed
 each $\mu$ in $X^+$ has a unique representative $\tilde\mu$ in $\mathbb Z^n$ with final entry 0. The cell
 corresponding to $\mu$ is $\mathbf c_{\tilde\mu}$, thought of as a cell of $\U$. (It is actually a cell of $\U$,
 $\mathbf U_D$, and $\mathfrak A_D$!)
\end{proof}

Note that this classification has an interesting consequence: The number of left cells in a two-sided cell
 $\mathbf c_\lambda$ of $\mathfrak A_D$ depends only on the element of $X^+$ it determines, as can be seen from
 the formula in Proposition \ref{leftcellnumber}. Thus since each left cell of $\mathfrak A_D$ intersects $\mathbf U_D$
 in at most one left cell, and the two algebras have the same number of left cells, this intersection is
 always nonempty.

We may also give an explicit formula for the value of the $a$ function, using the fact that we know the value
of the corresponding function on the Hecke algebra. Indeed if $w \in \hat A_{D-1}$ lies in the cell $\mathbf c_\lambda$
 then $a(w)=(D-\sum\lambda_i^2)/2$.

\begin{lemma}
Let $\mu \in X^+$ and let $\tilde\mu \in \mathbb Z^n$ be its representative with $0$ in the final entry.
 Suppose $b \in \B$ lies in the cell $\mathbf c_\mu$ corresponding to $\mu$, and $b1_\nu=b$, for some $\nu \in X$.
 Pick the unique representative $\upsilon$ of $\nu$ in $\mathbb Z^n$ such that $\sum_{i=1}^n \lambda_i = \sum_{i=1}^n \upsilon_i$.
 Then we have we have $a(b)=\sum_{i=1}^n(\lambda_i^2 -\upsilon_i^2)$.
\hfill $\square$
\end{lemma}

We have also already constructed the asymptotic algebra $A_{\mu}$ for each $\mu \in X^+$.
 This is just the ring $J_{\mathbf c_{\tilde\mu}}$ constructed in the previous section, where $\tilde\mu$ is the
 representative of $\mu$ described above.

Let $G_\mu \colon= \prod_{i=1}^{n-1} GL_{\mu(i)}(\mathbb C)$, and let $R_\mu$ be its representation ring.
 Combining the above with Proposition \ref{Jstructure} we find that the asymptotic ring $A_\mu$ is isomorphic to
 a matrix ring over $R_\mu$ of size $\prod_{i=1}^{n-1} \binom{n}{i}^{\mu(i)}$. Thus we may pull-back modules
 of this matrix ring to obtain modules for $\U$. These are the ``extremal weight modules'' of Kashiwara \cite{K02},
 which are in turn related to the universal standard modules defined by Nakjima in his geometric classification
 of simple modules for quantum affine algebras. Indeed Kashiwara has a number of conjectures about the structure of
 these modules which he suggests should be closely related to the conjectures of Lusztig that we establish here for
 $\widehat{\mathfrak{sl}_n}$ (see the remark below). Kashiwara's conjectures have recently been proved in the
 simply-laced case by Nakajima \cite{N} and Beck \cite{B}, and using them it is easy to show that the modules
 obtained from the asymptotic algebra are indeed the extremal weight modules. This observation can be used to give
 another proof of the formula for the number of left cells in a two--sided cell.
 It should be possible to give another approach to the results of this paper using these techniques, which would
 work for all the simply-laced cases.

\medskip

\begin{remark}
The results of this section establish (in the case of $\widehat{\mathfrak{sl}}_n$) all the conjectures
in \cite[section 5]{L95}. It should be noted that paragraph $5.4$ of that section contains a misprint.
Given $\lambda \in X^+$ the numbers $\lambda(i)$, for $i \in I_0$, should be given by the formula
$\lambda(i)= \langle i,\lambda \rangle$.
\end{remark}

\textit{Acknowledgements:} The results in this paper form part of my thesis written under the supervision
of George Lusztig. I would like to thank him for all his encouragement and advice, on the subject of this
paper and much else besides. This paper was written up while the author was supported by the Clay Mathematics
Institute.

\end{document}